# A Law of Likelihood for Composite Hypotheses

Zhiwei Zhang

*Food and Drug Administration, Rockville, Maryland, USA*

Summary

The law of likelihood underlies a general framework, known as the likelihood paradigm, for representing and interpreting statistical evidence. As stated, the law applies only to simple hypotheses, and there have been reservations about extending the law to composite hypotheses, despite their tremendous relevance in statistical applications. This paper proposes a generalization of the law of likelihood for composite hypotheses. The generalized law is developed in an axiomatic fashion, illustrated with real examples, and examined in an asymptotic analysis. Previous concerns about including composite hypotheses in the likelihood paradigm are discussed in light of the new developments. The generalized law of likelihood is compared with other likelihood-based methods and its practical implications are noted. Lastly, a discussion is given on how to use the generalized law to interpret published results of hypothesis tests as reduced data when the full data are not available.

## 1 Introduction

A major part of statistics is to interpret observed data as statistical evidence, often in response to questions like "what do the data say?". Yet there is no consensus among statisticians on what constitutes statistical evidence and how to measure its strength. Hypothesis tests and posterior probabilities are commonly used to interpret and communicate statistical evidence with regard to two competing hypotheses. The Neyman-Pearson theory for testing hypotheses was developed under a decision-theoretic framework that attempts to answer such questions as "what should I do?". The theory can lead to serious logical inconsistencies when used to address the problem of representing and interpreting statistical evidence (Royall, 1997, Chapter 2). The Bayesian approach, on the other hand, is more appropriate for questions of the form "what should I believe?". A posterior distribution, which incorporates prior information as well as the data, may not provide an objective representation of evidence in the observed data alone. A proper concept of evidence is missing from standard statistical theories.

The missing concept of evidence can be found in what Hacking (1965) termed the law of likelihood (LL):

> If one hypothesis, $H_1$, implies that a random variable $X$ takes the value $x$ with probability $f_1(x)$, while another hypothesis, $H_2$, implies that the probability is $f_2(x)$, then the observation $X = x$ is evidence supporting $H_1$ over $H_2$ if $f_1(x) > f_2(x)$, and the likelihood ratio, $f_1(x)/f_2(x)$, measures the strength of that evidence.

This point of view has led to a likelihood paradigm for interpreting statistical evidence, which carefully distinguishes evidence from error probabilities and personal belief (Royall, 1994, 1997; Blume, 2002). Royall (1997) also proposes benchmarks for the strength of statistical evidence. Specifically, a likelihood ratio (LR) exceeding $k = 8$ is considered fairly strong evidence, while $k = 32$ is used to define strong evidence. Royall (2000) analyzes the probability of observing misleading evidence under parametric models, and Blume (2008) provides a parallel analysis for sequential trials. Royall and Tsou (2003) and Blume et al. (2007) develop adjusted likelihood functions with certain robustness properties under model failure. Zhang (2008) advocates the use of empirical likelihood functions in nonparametric and semiparametric situations.

Most of the discussion so far about the likelihood paradigm has been limited to simple hypotheses. This is in sharp contrast to the tremendous relevance of composite hypotheses in statistical applications. In confirmatory clinical trials, for example, the primary objective is often to demonstrate that a new treatment is superior to a placebo/sham control, or not inferior to the standard of care by more than a specified amount. Also, there are bioequivalence trials designed to show that a generic drug is similar to a brand drug with respect to pharmacokinetic parameters. While Royall (1997, 2000) has considered some very special types of composite hypotheses (to be described in Section 2), he and Blume (2002) have not encouraged efforts to further extend the LL to general composite hypothese, due to concerns that will be addressed in Section 5. Blume (2002) argues that a graph of the likelihood function suffices for an evidential analysis and "no further reduction or summarization of the evidence is necessary". However, in many situations including the clinical studies mentioned earlier, it is important to not only look at the graph but also decide explicitly whether there is evidence supporting one specific hypothesis over another and, if so, how strong that evidence is. The LL makes this possible for simple hypotheses, and a reasonable generalization of the law for composite hypotheses will certainly be helpful. He et al. (2007) consider this problem in a finite parameter space without suggesting a practical solution. This article proposes a solution for general composite hypotheses.

The proposed solution will be derived in an axiomatic fashion in the next section, illustrated with real examples in Section 3, and analyzed asymptotically in Section 4. Section 5 discusses previous concerns about including composite hypotheses in the likelihood paradigm. Section 6 compares the proposed method with other likelihood-based methods. Lastly, a discussion is given in Section 7 on how to interpret published results of hypothesis tests as reduced data when the full data are not available.

## 2 Generalizing the law of likelihood

Let $X$ represent the data and suppose $X$ follows a distribution with density $f(\cdot; \theta)$, which is known up to a parameter (vector) $\theta$ taking values in $\Theta$. Then the likelihood for $\theta$, based on the observation $X = x$, is given by $L(\theta) = f(x; \theta)$. According to the LL, the data provide evidence supporting one parameter value $\theta_1$ over another value $\theta_2$ if $L(\theta_1) > L(\theta_2)$, and the strength of that evidence is measured by the LR $L(\theta_1)/L(\theta_2)$.

For this purpose, it is irrelevant whether the values $\theta_1$ and $\theta_2$ are predetermined or data-driven.

What do the data say about general hypotheses of the form $H_1 : \theta \in \Theta_1 \subset \Theta$ versus $H_2 : \theta \in \Theta_2 \subset \Theta$? Only some special cases have been considered. Obviously, if the likelihood $L$ is constant over each hypothesis with respective values $L_1$ and $L_2$, then the two hypotheses should be compared on the basis of the LR $L_1/L_2$ (Royall, 1997, Sections 1.7, 1.8). It has also been suggested that if the images $L(\Theta_1)$ and $L(\Theta_2)$ are intervals that do not overlap, then the evidence supports the hypothesis with larger likelihood values over the other one, with little discussion on how to measure the strength of that evidence (Royall, 2000; He et al., 2007). We shall take as an axiom a slight generalization of the latter suggestion.

**Axiom 1.** *If $\inf L(\Theta_1) > \sup L(\Theta_2)$, then there is evidence supporting $H_1$ over $H_2$.*

Also, it seems reasonable to expect evidential interpretations to be logically consistent, in the following sense.

**Axiom 2.** *If there is evidence supporting $H_1^*$ over $H_2$ and $H_1^*$ implies $H_1$, then the evidence also supports $H_1$ over $H_2$.*

These axioms together suggest the following generalization of the LL.

**Theorem 1.** *Let Axioms 1 and 2 hold. If $\sup L(\Theta_1) > \sup L(\Theta_2)$, then there is evidence supporting $H_1$ over $H_2$.*

*Proof.* Let $\sup L(\Theta_1) > \sup L(\Theta_2)$; then $\Theta_1$ is not empty and there exists $\theta_1 \in \Theta_1$ such that $L(\theta_1) > \sup L(\Theta_2)$. It follows from Axiom 1 that the simple hypothesis $H_1^* : \theta = \theta_1$ is supported over $H_2$. The result now follows directly from Axiom 2. □

It also seems natural to use the generalized likelihood ratio (GLR) $\sup L(\Theta_1)/\sup L(\Theta_2)$ to measure the strength of the evidence, although this does not follow from the axioms. This generalized law of likelihood (GLL) is consistent with the original law for simple hypotheses, and with previous suggestions for special cases of composite hyotheses (Royall, 1997, 2000). Moreover, the GLL coincides with the profile likelihood approach in the presence of nuisance parameters. Suppose $\theta = (\gamma, \omega)$, where $\gamma$ is of primary interest and $\omega$ is a nuisance parameter. While there may be ad hoc solutions available depending on the specific problem, a viable general approach is to represent the evidence about $\gamma$ with the profile likelihood $\widetilde{L}(\gamma) = \sup_\omega L(\gamma, \omega)$, which has good properties that justify its use (Royall, 2000, Section 5). The GLL, treating a "simple" hypothesis about $\gamma$ as a composite one about $(\gamma, \omega)$, would lead to the same answer for comparing two values of $\gamma$. Of course, unlike the profile likelihood approach, the GLL also applies to arbitrary composite hypotheses concerning $(\gamma, \omega)$.

As a by-product, the GLL allows one to assess the "absolute" evidence about a single parameter value or a single hypothesis, using its complement as the default comparator. This would not be possible without a mechanism to deal with composite hypotheses (Royall, 1997, 2000). Under the GLL, a hypothesis $H_1 : \theta \in \Theta_1$ is supported by the data if $\sup L(\Theta_1) > \sup L(\Theta_1^c)$, where the superscript $c$ denotes complement, and the strength of that evidence is measured by the GLR $\sup L(\Theta_1)/\sup L(\Theta_1^c)$. In particular, the evidence about a single parameter value, say $\theta_1$, is represented by the ratio $L(\theta_1)/\sup_{\theta \neq \theta_1} L(\theta)$. As we shall see in the examples to follow, excluding a single parameter value does not usually alter the supremum of the likelihood, which means $L(\theta_1)/\sup_{\theta \neq \theta_1} L(\theta) = L(\theta_1)/\sup L(\Theta) \leq 1$. Thus it is usually impossible to obtain empirical evidence supporting a single parameter value in a smooth model over its complement.

## 3 Real examples

The GLL will now be illustrated with some real examples. The first example, taken from Royall (1997, Section 1.9), is a clinical study in which 17 subjects were enrolled and given a new treatment. The outcome of interest was a binary indicator of treatment success, for which a binomial model is assumed. The success probability for the new treatment, denoted by $\theta$, was to be compared with the same probability for a standard treatment which was believed to be about 0.2. This gave rise to two composite hypotheses of interest: $H_1 : \theta \leq 0.2$ versus $H_2 : \theta > 0.2$. At the end of the study, nine subjects were found to have achieved treatment success and the resulting likelihood for $\theta$ is plotted in Figure 1. (In this and the subsequent likelihood plots, each likelihood function is divided by its maximum value so the peak value is invariably 1.) Royall (1997) recognizes that the LL does not allow us to compare $H_1$ and $H_2$ directly, and suggests that we should instead make pairwise comparisons between selected parameter values in $H_1$ and values in $H_2$. However, the choice of parameter values for pairwise comparisons can be rather subjective, and even after a large number of pairwise comparisons it may remain unclear how to answer the original question: Do the data support $H_1$ or $H_2$ and by how much? An application of the GLL yields a direct, unambiguous and objective answer: $H_2$ is supported over $H_1$ with a GLR of 91, which indicates strong evidence.

[Figure 1 about here.]

The second example is a randomized clinical study first reported by Rodary et al. (1989) and later discussed by many authors. The trial enrolled 164 children with nephroblastoma, who were randomly assigned to either chemotherapy or radiation therapy. The primary objective of the trial was to demonstrate that chemotherapy is non-inferior to radiation therapy with respect to the response rate. More precisely, non-inferiority here means that the response rate for chemotherapy is not lower than that for radiation therapy by more than a margin of 10%, which is considered the smallest clinically meaningful difference between two groups. The observed response rates were 94.3% (83/88) for chemotherapy and 90.8% (69/76) for radiation therapy. The profile likelihood function for the difference between the response rates in the two groups is computed using a bisection method as in Zhang (2006) and plotted in Figure 2. Under the GLL, the non-inferiority hypothesis is strongly supported with a GLR of 138. In fact, with a higher observed response rate in the chemotherapy group, there is even evidence supporting the superiority of chemotherapy to radiation therapy. This latter piece of evidence is rather weak, though, with a GLR of 1.4.

[Figure 2 about here.]

Finally, let us consider a bioequivalence trial described in Wellek (2003, Chapter 9). Bioequivalence trials are conducted to show that a generic drug or new formulation (test) is nearly equivalent in bioavailability to an approved brand-name drug or formulation (reference). There are different ways to measure bioavailability, but this example is primarily concerned with the area under the curve (AUC) for the serum concentration of the drug changing over time. The trial involved 25 patients and followed a crossover design where each patient was randomly assigned to a treatment sequence (test followed by reference or the opposite, with equal probabilities). Let $(Y_T, Y_R)$ denote the log-transformed AUC measurements in the test and reference periods, respectively, on the same subject. Following Choi et al. (2008), we assume that there are no sequence or period effects and that the measurements follow a simple bivariate normal model:

$$\begin{pmatrix} Y_T \\ Y_R \end{pmatrix} \sim N \left( \begin{pmatrix} \mu_T \\ \mu_R \end{pmatrix}, \begin{pmatrix} \sigma_T^2 & \rho \sigma_T \sigma_R \\ \rho \sigma_T \sigma_R & \sigma_R^2 \end{pmatrix} \right).$$

In this setting, it is natural to assess bioequivalence by comparing $\mu_T$ with $\mu_R$ and $\sigma_T$ with $\sigma_R$. Figures 3 and 4 display the profile likelihood functions for $\mu_T - \mu_R$ and $\sigma_T/\sigma_R$, respectively, both based on formulas given by Choi et al. (2008, Appendix A). In terms of the means, bioequivalence is usually defined as $|\mu_T - \mu_R| < 0.223$, which corresponds to $0.8 < \exp(\mu_T)/\exp(\mu_R) < 1.25$. As is clear in Figure 3, this bioequivalence hypothesis enjoys overwhelming support by the observed data, with a GLR greater than $10^6$. There is not a general definition of bioequivalence in terms of the standard deviations. One might, however, follow the same reasoning about the exponentiated means and consider the standard deviations close enough if $0.8 < \sigma_T/\sigma_R < 1.25$. The evidence regarding this latter hypothesis is largely neutral, with a GLR of 1.1.

[Figure 3 about here.]

[Figure 4 about here.]

## 4 Large-sample theory

In this section we consider the behavior of the GLL in large samples. Specifically, let $X = (Y_1, \ldots, Y_n)$, a collection of independent copies of some random variable $Y$, and let the density of $Y$ be modeled as $g(\cdot; \theta)$, $\theta \in \Theta$. Then the likelihood for $\theta$, based on the observations $Y_i = y_i$, $i = 1, \ldots, n$, is given by

$$L_n(\theta) = \prod_{i=1}^n g(y_i; \theta).$$

It is instructive to begin with the simple case of two simple hypotheses: $H_1 : \theta = \theta_1$ versus $H_2 : \theta = \theta_2$. The law of large numbers implies that, with probability 1,

$$l_n(\theta) := \frac{1}{n} \log L_n(\theta) = \frac{1}{n} \sum_{i=1}^n \log g(Y_i; \theta) \to E \log g(Y; \theta) =: l_\infty(\theta).$$

It follows that $L_n(\theta_1)/L_n(\theta_2)$ tends to $\infty$ if $l_\infty(\theta_1) > l_\infty(\theta_2)$ and to 0 if $l_\infty(\theta_1) < l_\infty(\theta_2)$. Thus the LL essentially orders the values in $\Theta$ according to $l_\infty$. If $l_\infty$ has a unique maximum at $\theta_0$, then $\theta_0$ will eventually dominate any other fixed value in $\Theta$. Obviously, $\theta_0$ is just the true value of $\theta$ if the model $g(\cdot;\theta)$ is correctly specified and suitably identified. Under an incorrect model, $\theta_0$ may be considered the "object of inference" (Royall and Tsou, 2003).

Under the GLL, the subsets of $\Theta$ are ordered according the suprema of their images under $l_\infty$. This can be formalized as follows.

**Theorem 2.** *Let $\Theta_1, \Theta_2$ be subsets of $\Theta$. If the collection of functions $\{\log g(Y;\theta) : \theta \in \Theta_j\}$ is Glivenko-Cantelli for $j = 1, 2$, then, with probability 1,*

$$\frac{\sup L_n(\Theta_1)}{\sup L_n(\Theta_2)} \to \begin{cases} \infty & \text{if } \sup l_\infty(\Theta_1) > \sup l_\infty(\Theta_2), \\ 0 & \text{if } \sup l_\infty(\Theta_1) < \sup l_\infty(\Theta_2). \end{cases}$$

This follows directly from the uniform law of large numbers. An extensive discussion of the Glivenko-Cantelli property, including techniques for its verification, can be found in van der Vaart and Wellner (1996). Certainly, any parameter set that contains $\theta_0$ (if it exists) attains the global maximum of $l_\infty$. On the other hand, a set that does not contain $\theta_0$ may or may not attain the global maximum, depending on certain properties of $l_\infty$. Some general characterizations are given below.

**Lemma 1.** *Suppose $l_\infty$ is continuous and maximized at $\theta_0$. If $\theta_0 \in \overline{\Theta}_1$ (closure of $\Theta_1$), then $\sup l_\infty(\Theta_1) = \max l_\infty(\Theta)$.*

**Lemma 2.** *Suppose $\theta_0$ is a well-separated maximizer of $l_\infty$ in the sense that for every $\epsilon > 0$,*

$$\sup_{\theta: \|\theta - \theta_0\| \geq \epsilon} l_\infty(\theta) < l_\infty(\theta_0).$$

*If $\theta_0 \notin \overline{\Theta}_1$, then $\sup l_\infty(\Theta_1) < \max l_\infty(\Theta)$.*

The preceding discussion leaves open the case that $\sup l_\infty(\Theta_1) = \sup l_\infty(\Theta_2)$, which can happen if $\theta_0$ lies on the boundary between $\Theta_1$ and $\Theta_2$ or, more generally, if $\theta_0 \in (\overline{\Theta}_1 \cap \overline{\Theta}_2)$. Take the first example in Section 3 comparing $H_1 : \theta \leq 0.2$ with $H_2 : \theta > 0.2$. If the true value of $\theta$ is 0.1 or 0.5, then the above results show that $H_1$ or $H_2$, respectively, will eventually dominate the other hypothesis. However, if $\theta = 0.2$, then the two hypotheses are tied with respect to $\sup l_\infty$ and a closer examination is required. The following theorem characterizes the asymptotic distribution of the GLR $\sup L_n(\Theta_1)/\sup L_n(\Theta_2)$ in terms of the limits of the sets $A_{jn} = \sqrt{n}(\Theta_j - \theta_0)$, $j = 1, 2$. A sequence of sets $A_n$ converges to a set $A$ if $A$ is the set of all limits $\lim a_n$ of convergent sequences $(a_n)$ with $a_n \in A_n$ for every $n$ and, moreover, the limit $a = \lim_k a_{n_k}$ of every convergent sequence $(a_{n_k})$ with $a_{n_k} \in A_{n_k}$ for every $k$ is contained in $A$. Also define the distance between a vector $b$ and a set $A$ in the same Euclidean space as $\|b - A\| = \inf_{a \in A} \|b - a\|$.

**Theorem 3.** *Suppose the model $\{g(\cdot;\theta) : \theta \in \Theta\}$ is differentialble in quadratic mean at $\theta_0$ in the sense of van der Vaart (1998, Section 7.2) with non-singular Fisher information matrix $I_{\theta_0}$. Suppose there exists a measurable function $h$ such that $E_{\theta_0}\{h(Y)^2\} < \infty$ and that, for every $\theta_1$ and $\theta_2$ in a neighborhood of $\theta_0$,*

$$|\log g(y;\theta_1) - \log g(y;\theta_2)| \leq h(y)\|\theta_1 - \theta_2\|.$$

Let $\Theta_1$, $\Theta_2$ be subsets of $\Theta$ for which the restricted maximum likelihood estimators $\widehat{\theta}_{jn} = \arg\max_{\Theta_j} L_n$, $j = 1, 2$, are consistent under $\theta_0$ and the sequence of sets $A_{jn}$ (defined above) converges to some $A_j$, $j = 1, 2$. Then, under $\theta_0$, we have

$$2\{\log\sup L_n(\Theta_1) - \log\sup L_n(\Theta_2)\} \xrightarrow{d} \|I_{\theta_0}^{1/2}W - I_{\theta_0}^{1/2}A_2\|^2 - \|I_{\theta_0}^{1/2}W - I_{\theta_0}^{1/2}A_1\|^2,$$

for $W$ normally distributed with mean 0 and variance matrix $I_{\theta_0}^{-1}$.

This result parallels Theorem 16.7 of van der Vaart (1998) for likelihood ratio tests and can be proved using similar arguments. Unlike Theorem 2, Theorem 3 does require correct specification of the model $g(\cdot; \theta)$. The differentiability condition for the model $g(\cdot; \theta)$ is satisfied for most models used in practice. The consistency of the restricted maximum likelihood estimators $\widehat{\theta}_{jn}$ typically requires that $\theta_0 \in (\Theta_1 \cap \Theta_2)$. In the more general case $\theta_0 \in (\overline{\Theta}_1 \cap \overline{\Theta}_2)$, one could replace each $\Theta_j$ with $\overline{\Theta}_j$ when applying the theorem and the conclusion would still hold for the original $\Theta_j$ as long as $L_n$ is continuous. To understand the limiting distribution given in Theorem 3, note that $I_{\theta_0}^{1/2}W$ is a standard normal vector of the same dimension as $\theta$. If $\theta_0$ is interior to $\Theta_j$, then $A_j$ is the entire Euclidean space and $I_{\theta_0}^{1/2}V - I_{\theta_0}^{1/2}A_j \equiv 0$. The situation is more complicated if $\theta_0$ lies on the boundary of a parameter set. Consider, again, the first example in Section 3. When $\theta_0 = 0.2$ in this example, $A_1$ and $A_2$ are the negative and positive halflines respectively, and the GLR $\sup L_n(\Theta_1)/\sup L_n(\Theta_2)$ converges in distribution to $\exp[\{\mathbb{I}(Z < 0) - \mathbb{I}(Z > 0)\}Z^2/2]$, where $Z$ is a standard normal variable and $\mathbb{I}(\cdot)$ is the indicator function. Thus $H_1$ and $H_2$ are virtually symmetric in the limit even though one is technically correct and the other is not. As a more extreme example, suppose $\theta_0$ is an interior point of $\Theta$ and take $\Theta_1 = \{\theta_0\}$ and $\Theta_2 = \Theta_1^c$. Then $A_1 = \{0\}$ and $A_2$ equals the entire Euclidean space, so the limit in Theorem 3 is $-\chi^2_{\dim(\theta)}$. This result adds to the discussion at the end of Section 2 with an asymptotic approximation.

## 5 Discussion of previous concerns

Both Royall (1997) and Blume (2002) have emphasized that the LL is silent about composite hypotheses. Neither author has expressed much optimism about expanding the likelihood paradigm to include composite hypotheses, despite their practical relevance. The main concerns expressed in Royall (1997, Section 1.9) and Blume (2002, Section 2.6) are outlined and addressed below.

### 5.1 Logical evidence versus statistical evidence

Royall (1997) points out the distinction between logical evidence and statistical evidence, and argues that the former should not be substituted for the latter. This can be illustrated with the following example from Royall (1997, Section 1.7.2). Suppose that $\Theta_1 = \{\theta_1\}$ and $\Theta_2 = \{\theta_1, \theta_2\}$ for distinct parameter values $\theta_1$ and $\theta_2$. Because $\Theta_1 \subset \Theta_2$, the hypothesis $H_2 : \theta \in \Theta_2$ certainly appears more plausible than $H_1 : \theta \in \Theta_1$. However, if $L(\theta_1) = L(\theta_2)$ then the LL says that $H_1$ and $H_2$ are equally well supported by the data. This is not surprising because the LL concerns the statistical evidence alone, irrespective of the logical relationship between the two hypotheses.

On the other hand, it is reasonable to require that statistical evidence be interpreted in a logically consistent manner. In the example of the above paragraph, it seems difficult to imagine a data configuration that can be naturally interpreted as supporting $H_1$ over $H_2$, even if $L(\theta_1)$ and $L(\theta_2)$ are allowed to differ in any possible way. A theory that allows such an illogical conclusion would be really troubling. This is the rationale for imposing Axiom 2 in developing the GLL.

The GLL respects logical relationships among hypotheses without substituting logical evidence for statistical evidence. This can be illustrated with the same example discussed in the above two paragraphs. Under the GLL, there can never be strict support for $H_1$ over $H_2$, eliminating logical inconsistencies. On the other hand, the GLL does not confuse the logical relationship $H_1 \Rightarrow H_2$ with statistical evidence. It never lends evidential support to $H_2$ over $H_1$ without a solid statistical basis, because the GLR of $H_2$ to $H_1$ cannot be greater than 1 unless $L(\theta_2) > L(\theta_1)$. This can also be seen in another example from Royall (1997, Section 1.9) with $\Theta_1 = \{\theta_1, \theta_3\}$, $\Theta_2 = \{\theta_2\}$ and $L(\theta_1) < L(\theta_2) < L(\theta_3)$. Based on the latter inequality, $H_1$ is supported over $H_2$ under the GLL. This conclusion will not hold if the likelihood becomes such that $L(\theta_1) < L(\theta_2) = L(\theta_3)$ even though the logical relationship between $H_1$ and $H_2$ is unchanged. More generally, for any hypothesis $H_1 : \theta \in \Theta_1$ to be supported over another hypothesis $H_2 : \theta \in \Theta_2$, the GLL requires the existence of $\theta_1 \in \Theta_1$ such that $L(\theta_1) > \sup L(\Theta_1)$, a crucial piece of statistical evidence that cannot be replaced by any logical evidence.

## 5.2 The role of a prior distribution

If a prior distribution can be specified for $\theta$, then the probabilities $P(X = x|H_1)$ and $P(X = x|H_2)$ can be evaluated and the LL can be used to assess the evidence about $H_1$ versus $H_2$. As Royall (1997) and Blume (2002) point out, this Bayesian approach reduces composite hypotheses into simple ones by modifying the probability model. The choice of a prior distribution can be arbitrary, and an LR that involves a prior distribution may not provide an objective representation of the observed evidence.

Given his insight into the Bayesian approach, it is interesting that Royall supports his reservation about considering composite hypotheses in the likelihood paradigm with the following Bayesian observation. He observes that the posterior probability ratio $P(H_1|X = x)/P(H_2|X = x)$ can be larger or smaller than the prior probability ratio $P(H_1)/P(H_2)$, depending on the prior distribution, unless one hypothesis dominates the other as in Axiom 1 (Royall, 1997, Section 1.9). If anything, this appears to highlight the subjectivity inherent in the Bayesian approach. In no way does the above observation suggest that statistical evidence about composite hypotheses cannot be interpreted objectively.

## 5.3 The lack of a unique solution

Blume's main concern about composite hypotheses appears to be that there is no unique way to deal with them (Blume, 2002, Section 2.6). This, of course, is the case with many problems, including the very problem the LL aims to address. (Statisticians and scientists who have not yet accepted the LL wholeheartedly need not regard it as the unique approach to statistical evidence.) Within the likelihood paradigm, there are different ways to obtain a likelihood function for the parameter of interest in the presence of nuisance parameters (a special type of composite hypotheses). Among other possibilities, the profile likelihood approach has been shown to have desirable statistical properties and appears to be a viable general approach to dealing with nuisance parameters (Royall, 1997, 2000; Blume, 2002). The GLL provides a natural extension of the profile likelihood approach to more general composite hypotheses.

It is actually desirable to have a set of possible solutions to choose from. Besides maximizing the likelihood over each hypothesis as in this paper, Blume (2002) also mentions other possible approaches to composite hypotheses (to illustrate the lack of a unique solution). One of the alternatives mentioned is to minimize the likelihood over each hypothesis, which clearly makes no sense. Suppose, for example, that $L(\theta_1) > 0$ for some $\theta_1$ and that $\inf L(\Theta) = 0$, which is not uncommon. The minimum likelihood rule would then indicate (infinitely) strong evidence supporting $H_1 : \theta = \theta_1$ over $H_2 : \theta \in \Theta$, even though the latter hypothesis is trivially true. Blume also mentions the possibility of averaging the likelihood over a composite hypothesis with a weight function, which is essentially the Bayesian approach. As noted earlier, the Bayesian approach relies on external information and may not provide an objective representation of the observed evidence. In contrast, The GLL provides an objective representation of evidence and avoids illogical and counterintuitive conclusions.

# 6 Comparison with other procedures

## 6.1 Connection with likelihood ratio tests

As we have seen in the examples discussed so far, $\Theta_1$ and $\Theta_2$ need not be disjoint or exhaustive for the GLL to apply. There are, however, many applications where it is customary to take $\Theta_2 = \Theta_1^c$. These problems are usually tackled with statistical tests, such as likelihood ratio tests (LRTs). To be specific, let $H_1 : \theta \in \Theta_1$ be considered the null hypothesis and $H_2 : \theta \notin \Theta_1$ the alternative; this suggests that evidence supporting $H_2$ over $H_1$ is of particular interest. The LRT is based on the statistic

$$\frac{\sup L(\Theta)}{\sup L(\Theta_1)} = \frac{\sup L(\Theta_1) \vee \sup L(\Theta_2)}{\sup L(\Theta_1)} = \frac{\sup L(\Theta_2)}{\sup L(\Theta_1)} \vee 1,$$

where $\vee$ denotes maximum. The null hypothesis will formally be rejected if the above test statistic is greater than some critical value, which is typically greater than 1. Thus, for the purpose of seeking evidence for $H_2$, the LRT statistic is essentially equivalent to the GLR $\sup L(\Theta_2)/\sup L(\Theta_1)$ given by the GLL.

However, the same GLR can be interpreted in different ways. In the likelihood paradigm, the GLR is all that is needed to compare the two hypotheses. It determines both the nature and the strength of the evidence, and can be placed on a universal scale together with GLRs from different problems. The strength of statistical evidence is measured on a continuum and should be understood as such, even though descriptive benchmarks are sometimes used to facilitate communication. In an LRT, the test statistic is to be compared with a critical value before a dichotomous conclusion (whether to reject $H_1$ in favor of $H_2$) can be reached. The critical value is derived from the (asymptotic) null distributions of the test statistic and usually depends on the significance level and certain properties of $\Theta_1$ (van der Vaart, 1998, Chapter 16). The use of a $p$-value provides some continuity and eliminates the dependence on the significance level, but a $p$-value still depends on other features of the problem that are irrelevant from an evidential point of view. In general, Royall (1997) argues that hypothesis tests are not appropriate tools for interpreting data as evidence.

## 6.2 Confidence sets versus support sets

Denote by $c_\alpha(\theta_1)$ the LRT critical value for testing $H_1 : \theta = \theta_1$ against $H_2 : \theta \neq \theta_1$ at significance level $\alpha$. The associated $1 - \alpha$ confidence set for $\theta$ is given by $\{\theta \in \Theta : L(\theta) > \sup L(\Theta)/c_\alpha(\theta)\}$. If $c_\alpha(\theta) \equiv c_\alpha$, as is often the case, then the confidence set is simply $\{\theta : L(\theta) > \sup L(\Theta)/c_\alpha\}$. This, interestingly, coincides with a likelihood support set. Royall (1997) and Blume (2002) have discussed likelihood support sets in one-dimensional situations, where the support sets are typically intervals. In general, the $1/k$ support set for $\theta$ can be defined as

$$S_k = \{\theta : L(\theta) > \sup L(\Theta)/k\}, \qquad k > 1. \tag{1}$$

Despite a similar appearance, a support set is to be interpreted differently than a confidence set. The usual interpretation of confidence sets in terms of long-run coverage does not fit well into the likelihood paradigm, where the emphasis is placed on understanding the observed data (as opposed to fictitious repetitions of the same experiment). The LL leads to the following interpretation of support sets. The values in the $1/k$ support set $S_k$ are "consistent with the observations" in the sense that no other value can be better supported by a factor greater than $k$ (Royall, 1997, Section 1.12). An alternative, perhaps more straightforward interpretation is made available by the GLL. Recall from Section 2 that evidence about a single hypothesis could be evaluated using its complement as the default comparator. In this sense, $S_k$ is simply the smallest parameter set supported by a factor of $k$.

**Theorem 4.** *If* $\sup L(S)/\sup L(S^c) \geq k$ *for some* $S \subset \Theta$, *then* $S_k \subset S$.

*Proof.* By assumption, $\sup L(S^c) \leq \sup L(S)/k \leq \sup L(\Theta)/k$. It follows that

$$S^c \subset \{\theta : L(\theta) \leq \sup L(\Theta)/k\} = S_k^c,$$

from which the result is immediate. □

## 6.3 Drawing support from support sets

Choi et al. (2008) propose to use profile likelihood support intervals to assess bioequivalence. As discussed in Section 3, the parameter of interest in this context is usually a scalar parameter $\gamma \in \Gamma$ that summarizes the difference in bioavailability between a test drug or formulation and a reference. Suppose bioequivalence is defined as $\gamma \in (\gamma_L, \gamma_U) = \Gamma_{BE}$ for specified bounds $\gamma_L$ and $\gamma_U$. Write $\theta = (\gamma, \omega)$ for the entire parameter vector that determines the distribution of the data, with $\omega \in \Omega$ considered the nuisance parameter. Choi et al. work with the profile likelihood $\widetilde{L}(\gamma) = \sup_\omega L(\gamma, \omega)$, from which support intervals for $\gamma$ can be derived as

$$\widetilde{S}_k = \{\gamma \in \Gamma : \widetilde{L}(\gamma) > \sup \widetilde{L}(\Gamma)/k\}, \qquad k > 1.$$

Note that $\sup \widetilde{L}(\Gamma) = \sup L(\Theta)$ and that $\widetilde{S}_k = \{\gamma : (\gamma, \omega) \in S_k \text{ for some } \omega\}$ with $S_k$ defined by (1). Choi et al. suggest that if $\widetilde{S}_k \subset \Gamma_{BE}$ then there is evidence supporting the bioequivalence hypothesis; the larger $k$ is, the stronger the evidence. They also suggest what is essentially a sensitivity analysis using several support intervals (say with $k = 5, 8, 32$). Implicit in the latter suggestion is an attempt to quantify the strength of the evidence with the "largest" $k$ for which the bioequivalence hypothesis is supported.

The GLL provides a general theoretical basis for the above suggestions. To see this, note first that $\widetilde{S}_k \subset \Gamma_{BE}$ if and only if $S_k \subset \Gamma_{BE} \times \Omega =: \Theta_{BE}$. Thus Choi et al.'s approach can be cast in terms of the parameter $\theta$ and the associated support sets $S_k$. The following theorem justifies their usage of support sets for an arbitrary hypothesis concerning $\theta$.

**Theorem 5.** *Let $\Theta_A \subset \Theta$ and write $r_A = \sup L(\Theta_A)/\sup L(\Theta_A^c)$. (a) If $S_k \subset \Theta_A$ for some $k > 1$, then $r_A \geq k$. (b) Denote $k^* = \sup\{k > 1 : S_k \subset \Theta_A\}$, which will be set to 1 if there is no qualifying $k$. Then $r_A > 1$ if and only if $k^* > 1$, in which case $r_A = k^*$.*

*Proof.* If $S_k \subset \Theta_A$, then
$$r_A = \frac{\sup L(\Theta_A)}{\sup L(\Theta_A^c)} \geq \frac{\sup L(S_k)}{\sup L(S_k^c)} \geq k,$$
proving statement (a). The "if" part of statement (b) follows directly from statement (a), which, by a limiting argument, further implies that $r_A \geq k^*$. To prove the "only if" part, we can invoke Theorem 4 with $S = \Theta_A$ and $k = r_A$, which also shows that $r_A \leq k^*$. The proof is complete upon combining the two arguments. □

## 7 Hypothesis tests as reduced data

When reporting their research findings, scientists do not always provide the raw data and sometimes only present the final results of statistical tests concerning their research hypotheses. Without access to the raw data, an interested reader may not be able to produce the likelihood function for the parameter of interest and is often forced to work with the results of hypothesis tests. Such difficulties do not necessarily force us out of the likelihood paradigm, as the GLL can still be used to interpret hypothesis tests as reduced data. This will be illustrated below for both dichotomous test results and $p$-values.

## 7.1 Dichotomous test results

Suppose the null hypothesis $H_1$ is to be tested against the alternative $H_2$ at significance level $\alpha$. Let $T = 1$ if $H_1$ is rejected in favor of $H_2$; 0 otherwise. The GLR of $H_2$ to $H_1$ based on the result $T = t$ is given by

$$r_T(t) = \sup_{H_2} P(T = t) / \sup_{H_1} P(T = t), \qquad t = 0, 1.$$

In the case of simple hypotheses, this has been discussed by Royall (1997, pp. 48–49). For an arbitrary $H_1$, $\sup_{H_1} P(T = 1)$ is the size of the test, which must not exceed $\alpha$ for a level-$\alpha$ test. In fact, except in very discrete cases, a reasonable test should have size close to $\alpha$. On the other hand, $\sup_{H_2} P(T = 1)$ is the maximum power of the test, which is generally greater than $\alpha$ and frequently equal to 1. Thus it seems reasonable to expect $r_T(1) > 1$, which justifies the usual interpretation of a rejection as evidence for the alternative. The GLR based on a non-rejection appears less predictable, consistent with the conventional wisdom that a failure to reject the null does not necessarily support the null. It should be noted that while the GLL appears consistent with conventional interpretations of hypothesis tests, it does not justify the use of such tests to reduce the data. The GLL should ideally be applied to the original data.

More can be said about the GLR $r_T$ for some typical hypotheses concerning a scalar $\theta \in (\underline{\theta}, \overline{\theta})$. Consider for instance the one-sided hypotheses $H_1 : \theta \leq \theta^*$ versus $H_2 : \theta > \theta^*$. In terms of the power function $\pi(\theta) = P_\theta(T = 1)$, the GLR can now be written as

$$r_T(t) = \begin{cases} \{1 - \inf_{\theta > \theta^*} \pi(\theta)\} / \{1 - \inf_{\theta \leq \theta^*} \pi(\theta)\}, & t = 0; \\ \sup_{\theta > \theta^*} \pi(\theta) / \sup_{\theta \leq \theta^*} \pi(\theta), & t = 1. \end{cases}$$

In many examples, $\pi(\theta)$ is increasing in $\theta$ with $\pi(\underline{\theta}+) = 0$, $\pi(\theta^*) = \alpha$ and $\pi(\overline{\theta}-) = 1$. If this is the case, then

$$r_T(t) = \begin{cases} 1 - \alpha, & t = 0; \\ 1/\alpha, & t = 1. \end{cases}$$

Thus, if $H_1$ is rejected, the resulting evidence supporting $H_2$ over $H_1$ is moderate ($r_T = 20$) for $\alpha = 0.05$ and strong ($r_T = 40$) for $\alpha = 0.025$. If $H_1$ is not rejected, then $H_1$ is supported over $H_2$ but the evidence is very weak for common values of $\alpha$. Sometimes the one-sided hypotheses are formulated as $H_1' : \theta = \theta^*$ versus $H_2 : \theta > \theta^*$. Such a reformulation does not usually require a different test. Assuming the same properties of $\pi$ as stated above, we then have

$$r_T(t) = \begin{cases} \{1 - \inf_{\theta > \theta^*} \pi(\theta)\} / \{1 - \pi(\theta^*)\} = 1, & t = 0; \\ \sup_{\theta > \theta^*} \pi(\theta) / \pi(\theta^*) = 1/\alpha, & t = 1. \end{cases}$$

The reformulated null hypothesis cannot be supported even if it is not rejected.

Suppose the simple null $H_1 : \theta = \theta^*$ is to be tested against the two-sided alternative $H_2 : \theta \neq \theta^*$. Assume the power function $\pi(\theta)$ equals $\alpha$ at $\theta = \theta^*$, increases as $\theta$ moves away from $\theta^*$, and tends to 1 as $\theta$ approaches $\underline{\theta}$ or $\overline{\theta}$. Then

$$r_T(t) = \begin{cases} 1, & t = 0; \\ 1/\alpha, & t = 1. \end{cases}$$

Again, the simple null is never supported.

Now let us consider an equivalence testing problem with hypotheses $H_1 : |\theta - \theta^*| \geq \delta$ and $H_2 : |\theta - \theta^*| < \delta$ for some $\delta > 0$. Here it may be reasonable to assume that $\pi(\theta)$ attains its maximum $\pi_{\max}$ at $\theta = \theta^*$, decreases as $\theta$ moves away from $\theta^*$, equals $\alpha$ at $\theta = \theta^* \pm \delta$, and tends to 0 as $\theta$ approaches $\underline{\theta}$ or $\overline{\theta}$. The maximum power $\pi_{\max}$ is usually less than 1 but greater than $\alpha$. In this case the GLR of $H_2$ to $H_1$ is given by

$$r_T(t) = \begin{cases} 1 - \alpha, & t = 0; \\ \pi_{\max}/\alpha, & t = 1. \end{cases}$$

This provides a realistic example where the maximum power of the test is relevant even after rejecting the null.

### 7.2  $p$-values

Let $U$ denote a $p$-value for testing the null hypothesis $H_1$ against the alternative $H_2$. Depending on the nature of $X$ and the procedure, $U$ may be discrete or continuous. In either case we write $f_U$ for the probability density of $U$ with respect to an appropriate measure. The GLR of $H_2$ to $H_1$ based on the result $U = u$ is

$$r_U(u) = \sup_{H_2} f_U(u) / \sup_{H_1} f_U(u), \qquad 0 < u < 1.$$

In practice, $U$ is often determined by a test statistic $V$ through a smooth monotone function, such as 1 minus a reference distribution function, in which case $r_U$ is equivalent to the analogous GLR based on $V$ (denoted by $r_V$).

For example, suppose $X = (Y_1, \ldots, Y_n)$ is a random sample from $N(\mu, \sigma^2)$ with $\sigma^2$ known, and consider the hypotheses $H_1 : \mu \leq 0$ versus $H_2 : \mu > 0$. In this situation it is common to use the $p$-value $U = 1 - \Phi(V)$, where $\Phi$ is the standard normal distribution function and $V = n^{-1/2} \sum_{i=1}^n Y_i / \sigma$. The corresponding GLR can be obtained as

$$r_U(u) = r_V(\Phi^{-1}(1-u)) = \frac{\sup_{\mu > 0} \phi(\Phi^{-1}(1-u) - \sqrt{n}\mu/\sigma)}{\sup_{\mu \leq 0} \phi(\Phi^{-1}(1-u) - \sqrt{n}\mu/\sigma)}$$

$$= \begin{cases} \exp(\{\Phi^{-1}(1-u)\}^2/2), & u \leq 0.5; \\ \exp(-\{\Phi^{-1}(1-u)\}^2/2), & u > 0.5, \end{cases} \qquad (2)$$

where $\phi$ denotes the standard normal density function. Note that $H_2$ is supported over $H_1$ if and only if $U < 0.5$. In fact, since $V$ is sufficient for $\mu$ in this example, the above GLR is the same as that based on the full data $X$. A less trivial example would be a two-sample problem with a random sample $(Y_{11}, \ldots, Y_{1n_1})$ from $N(\mu_1, \sigma^2)$ and another sample $(Y_{21}, \ldots, Y_{2n_2})$ from $N(\mu_2, \sigma^2)$. Assume $\sigma^2$ is known and consider the hypotheses $H_1 : \mu_2 \leq \mu_1$ versus $H_2 : \mu_2 > \mu_1$. Then the test statistic $V = (n_1^{-1} + n_2^{-1})^{-1/2}(\overline{Y}_2 - \overline{Y}_1)/\sigma$, with $\overline{Y}_j = n_j^{-1}\sum_{i=1}^{n_j} Y_{ij}$ ($j=1,2$), would not be sufficient. Nonetheless, the GLR of $H_2$ to $H_1$ based on the $p$-value $U = 1 - \Phi(V)$ continues to follow expression (2).

When interpreting a $p$-value, one might be tempted to invoke the characterization of $U$ as the "lowest" significance level at which $H_1$ is rejected, that is, $U = \inf\{\alpha : T_\alpha = 1\}$, where $T$ is subscripted to emphasize its dependence on $\alpha$. The observation $U = u$ implies that the use of any $\alpha > u$ would lead to a rejection and hence a GLR $r_{T_\alpha}(1)$, and it might seem natural to use the quantity $r_+(u) = \sup_{\alpha > u} r_{T_\alpha}(1)$ to represent the evidence in the observed $p$-value. While this quantity may be easy to compute when a simple expression for $r_{T_\alpha}(1)$ is available as in Section 7.1, its interpretation can be problematic. First, a symmetric argument based on the dual identity $U = \sup\{\alpha : T_\alpha = 0\}$ would lead to the dual quantity $r_-(u) = \inf_{\alpha < u} r_{T_\alpha}(0)$, which is generally different and typically smaller than $r_+(u)$. It is not clear how to reconcile the difference. More importantly, the GLRs $r_{T_\alpha}(t)$ are defined in Section 7.1 for a fixed $\alpha$. If $\alpha$ is allowed to depend on $U$, a random variable, then $T_\alpha$ will become a different statistic to which the discussion of Section 7.1 no longer applies.

# References

Blume, J. D. (2002). Likelihood methods for measuring statistical evidence. *Statistics in Medicine*, 21, 2563–2599.

Blume, J. D. (2008). How often likelihood ratios are misleading in sequential trials. *Communications in Statistics: Theory and Methods*, 37, 1193-1206.

Blume, J. D., Su, L., Olveda, R. M. and McGarvey, S. T. (2007). Statistical evidence for GLM regression parameters: A robust likelihood approach. *Statistics in Medicine*, 26, 2919–2936.

Choi, L., Caffo, B. and Rohde, C. (2008). A survey of the likelihood approach to bioequivalence trials. *Statistics in Medicine*, 27, 4874–4894.

Hacking, I. (1965). *Logic of Statistical Inference*. Cambridge University Press, New York.

He, Y., Huang, W. and Liang, H. (2007). Axiomatic development of profile likelihoods as the strength of evidence for composite hypotheses. *Communications in Statistics: Theory and Methods*, 36, 2695–2706.

Rodary, C., Com-Nougue, C. and Tournade, M. (1989). How to establish equivalence between treatments: A one-sided trial in paediatric oncology. *Statistics in Medicine*, 8, 593–598.

Royall, R. (1994). The elusive concept of statistical evidence. *Bayesian Statistics 4*, eds. J. M. Bernardo, J. O. Berger, A. P. Dawid, and A. F. M. Smith. Clarendon Press, Oxford, U.K.

Royall, R. (1997). *Statistical Evidence: A Likelihood Paradigm*. Chapman & Hall, Boca Raton, FL.

Royall, R. (2000). On the probability of observing misleading statistical evidence. *Journal of the American Statistical Association*, 95, 760–768.

Royall, R. and Tsou, T.-S. (2003). Interpreting statistical evidence by using imperfect models: Robust adjusted likelihood functions. *Journal of the Royal Statistical Society, Series B*, 65, 391–404.

van der Vaart, A. W. (1998). *Asymptotic Statistics*. Cambridge University Press, Cambridge.

van der Vaart, A. W. and Wellner, J. A. (1996). *Weak Convergence and Empirical Processes with Applications to Statistics*. Springer-Verlag, New York.

Wellek, S. (2003). *Testing Statistical Hypotheses of Equivalence*. CRC Press, Boca Raton, FL.

Zhang, Z. (2006). Non-inferiority testing with a variable margin. *Biometrical Journal*, 48, 948–965.

Zhang, Z. (2008). Interpreting statistical evidence with empirical likelihood functions. Under review.

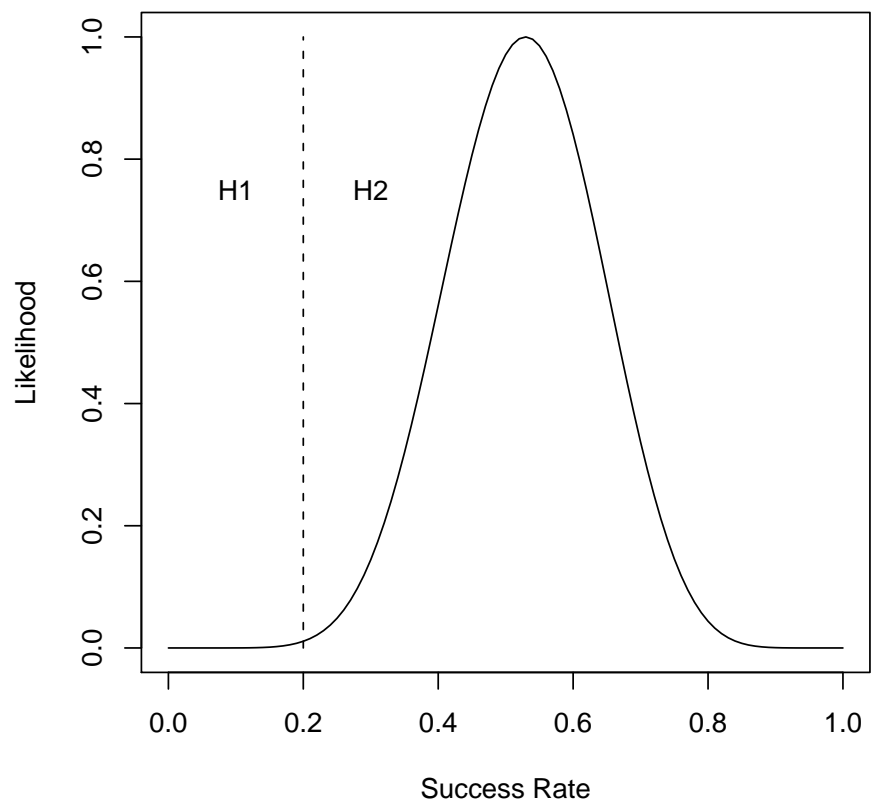

Figure 1: Likelihood function for the success rate $\theta$ in the first example of Section 3, with the dashed line separating the hypotheses $H_1 : \theta \leq 0.2$ and $H_2 : \theta > 0.2$.

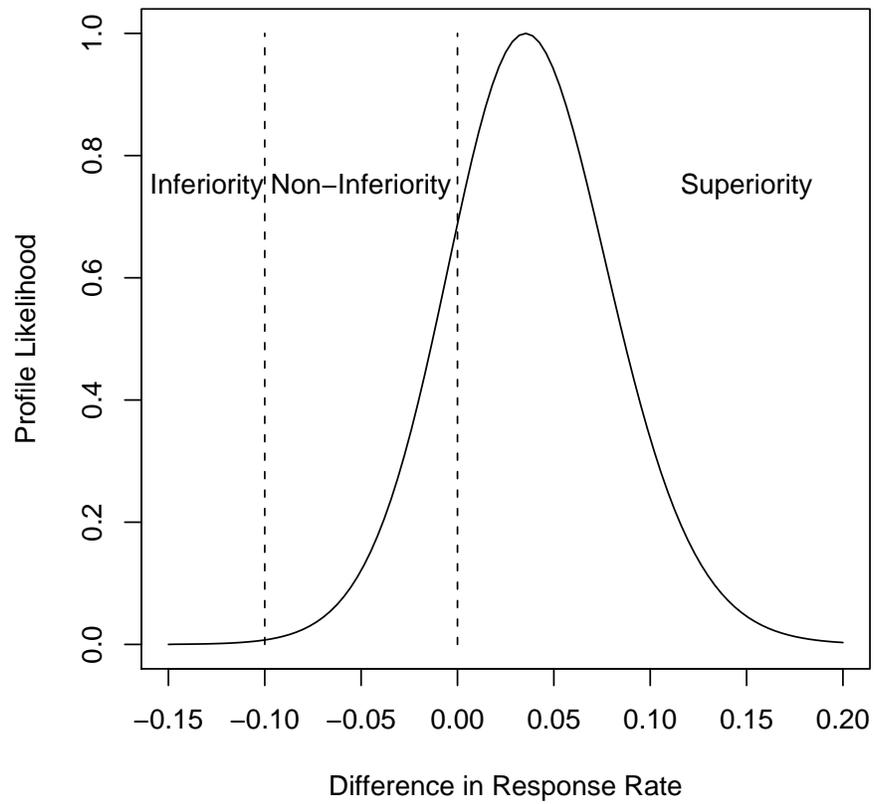

Figure 2: Profile likelihood function for the difference in response rate (chemotherapy − radiation therapy) in the second example of Section 3, with the dashed lines separating regions of inferiority, non-inferiority and superiority.

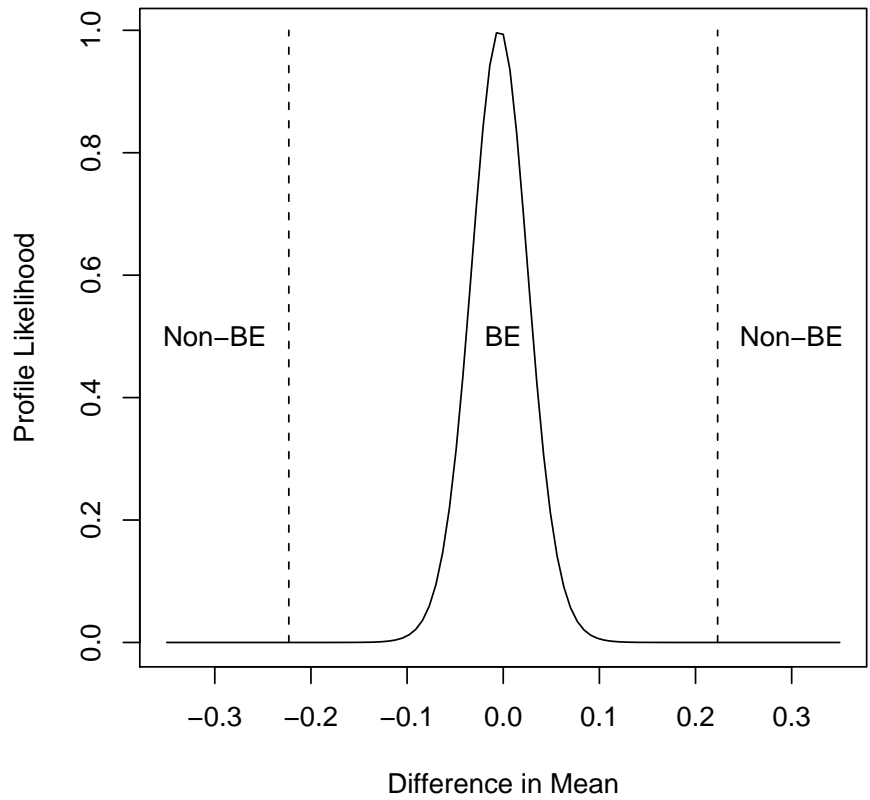

Figure 3: Profile likelihood function for the difference in mean log-AUC (test−reference) in the third example of Section 3, with the dashed lines separating regions of bioequivalence (BE) and non-BE.

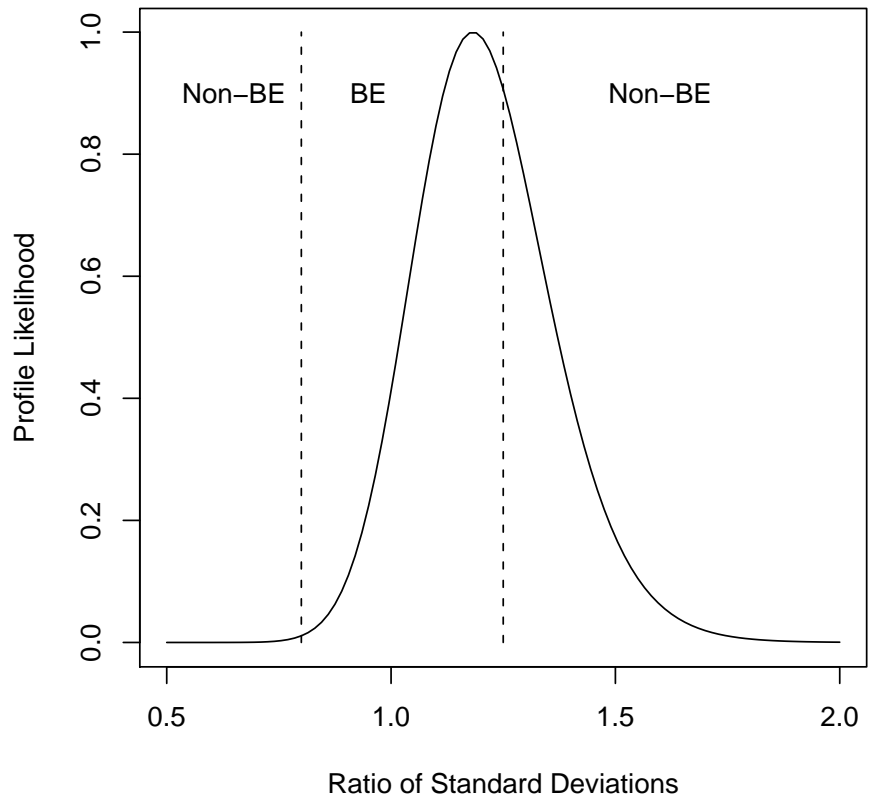

Figure 4: Profile likelihood function for the ratio of log-AUC standard deviations (test/reference) in the third example of Section 3, with the dashed lines separating regions of bioequivalence (BE) and non-BE.